\documentclass[11pt]{article}

\usepackage{amssymb,latexsym}
\usepackage{amsmath}
\usepackage{citesort}
\usepackage{eucal}

\title{The characteristic classes of Morita equivalent star products
  on symplectic manifolds} 

\author{\textbf{Henrique Bursztyn\thanks{henrique@math.berkeley.edu}}
  \\[0.5cm]
  Department of Mathematics\\
  UC Berkeley\\
  94720 Berkeley, CA, USA
  \\[1cm]
  {\bf Stefan Waldmann\thanks{Stefan.Waldmann@physik.uni-freiburg.de}} 
  \\[0.5cm]
  Fakult{\"a}t f{\"u}r Physik\\ 
  Albert-Ludwigs-Universit{\"a}t Freiburg\\
  Hermann Herder Stra{\ss}e 3\\
  D 79104 Freiburg\\
  Germany
  }

\date{August 2001}

\newcommand{\im} {\mathrm{i}}
\newcommand{\eu} {\mathrm{e}}

\newcommand{\cc}[1]      {\overline{{#1}}}
\newcommand{\supp}       {\mathop{{\mathrm{supp}}}}
\newcommand{\id}         {{\mathsf{id}}}
\newcommand{\ad}         {\mathop{{\mathrm{ad}}}}
\newcommand{\End}        {\mathop{{\mathsf{End}}}}
\newcommand{\ring}[1]    {{\mathsf{{#1}}}}
\newcommand{\SP} [1]     {{\left\langle{{#1}}\right\rangle}}
\newcommand{\Unit}       {\mathsf{1}}
\newcommand{\Def}        {{\mathsf{Def}}}

\newcommand{\Exp}        {\operatorname{\mathrm Exp}}
\newcommand{\Ln}         {\operatorname{\mathrm Ln}}

\newcommand{\qe}         {\boldsymbol{e}}
\newcommand{\qs}         {\boldsymbol{s}}
\newcommand{\qpsi}       {\boldsymbol{\Psi}}
\newcommand{\qphi}       {\boldsymbol{\Phi}}
\newcommand{\qh}         {\boldsymbol{h}}
\newcommand{\qH}         {\boldsymbol{H}}
\newcommand{\qU}         {\boldsymbol{U}}
\newcommand{\qV}         {\boldsymbol{V}}
\newcommand{\qA}         {\boldsymbol{A}}
\newcommand{\qB}         {\boldsymbol{B}}
\newcommand{\qt}         {\boldsymbol{t}}

\newcommand{\CH}         {\check{\mathrm{H}}}
\newcommand{\Pic}        {\mathsf{Pic}}
\newcommand{\HdR}        {\mathrm{H}_{\scriptscriptstyle \mathrm{dR}}}

\newcommand{\FDiff}      {\mathsf{F}}
\newcommand{\SymD}       {\mathsf{D}}

\newcommand{\TinyW}      {\mathrm{\scriptscriptstyle W}}

\newcommand{\starw}      {\mathbin{\star_{\TinyW}}}
\newcommand{\wrep}       {\varrho_{\TinyW}}
\newcommand{\TinyS}      {\mathrm{\scriptscriptstyle S}}
\newcommand{\stars}      {\mathbin{\star_{\TinyS}}}
\newcommand{\srep}       {\varrho_{\TinyS}}
\newcommand{\stark}      {\mathbin{\star_\kappa}}
\newcommand{\bulletw}    {\mathbin{\bullet_{\TinyW}}}
\newcommand{\bullets}    {\mathbin{\bullet_{\TinyS}}}

\newcommand{\sRep}       {{}^*\textrm{-}\mathsf{Rep}}

\newtheorem{lemma} {Lemma} [section]
\newtheorem{proposition} [lemma] {Proposition}
\newtheorem{theorem} [lemma] {Theorem}
\newtheorem{corollary} [lemma] {Corollary}
\newtheorem{definition}[lemma] {Definition}

\newtheorem{remark}[lemma]{Remark}

\newenvironment{proof}{{\sc Proof:}}{{\hspace*{\fill} $\square$\\}}

\numberwithin{equation}{section}

\begin{document}

\maketitle

\begin{abstract}
     In this paper we give a complete characterization of Morita equivalent
     star products on symplectic manifolds in terms of their characteristic
     classes: two star products $\star$ and $\star'$
     on $(M,\omega)$ are Morita equivalent if and only if 
     there exists a symplectomorphism $\psi:M\longrightarrow M$
     such that the
     relative class $t(\star,\psi^*(\star'))$ is $2 \pi \im$-integral.   
     For star products on cotangent bundles, we show that this integrality condition
     is related to Dirac's quantization condition for magnetic charges.
\end{abstract}

%
%

\section{Introduction}
\label{sec:intro}

The concept of Morita equivalence has played an important role in
different areas of mathematics (see \cite{landsman:2000:pre} for an
overview) since its introduction in the study of unital rings
\cite{morita:1958}. In applications of noncommutative geometry to
$M$-theory \cite{connes.douglas.schwarz:1998}, Morita equivalence was
shown to be related to physical duality \cite{schwarz:1998},
motivating the study of Morita equivalence of quantum tori
\cite{rieffel.schwarz:1999}. In this setting, the problem is to
characterize constant Poisson structures $\theta$ on the $n$-torus
$T^n$ that, after strict deformation quantization \cite{rieffel:1989},
give rise to Morita equivalent $C^*$-algebras $T_\theta$.

In this paper we address the problem of characterizing Morita
equivalent  algebras obtained from formal deformation
quantization of Poisson manifolds \cite{bayen.et.al:1978}
(see \cite{gutt:2000,sternheimer:1998:pre,weinstein:1994} for
surveys). In this approach to quantization, quantum algebras of
observables are defined by formal associative deformations (in the
sense of \cite{gerstenhaber:1964}) of classical Poisson algebras
known as \emph{star products}.

The problem of classifying Morita equivalent star products on a
Poisson manifold $(M,\pi_0)$ can be phrased in terms of a canonical
action $\Phi$ of the Picard group $\Pic(M) \cong H^2(M,\mathbb Z)$
on $\Def(M,\pi_0)$, the moduli space of equivalence classes of
differential star products on $(M,\pi_0)$ \cite{bursztyn:2001:pre}. 
The action $\Phi$ is defined by deformation quantization of line
bundles on $M$ \cite{bursztyn.waldmann:2000b}, and two star products
$\star$, $\star'$ are Morita equivalent 
(as unital $\mathbb{C}[[\lambda]]$-algebras)
if and only if there exists a Poisson diffeomorphism
$\psi:M\longrightarrow M$ such that the classes 
$[\star]$ and $[\psi^*(\star')]$ lie in the
same $\Phi$-orbit. The semiclassical limit of this action was
described in \cite[Thm.~5.11]{bursztyn:2001:pre}.

Let $(M,\omega)$ be a symplectic manifold. The main result of this
paper is that, under the usual identification
\cite{nest.tsygan:1995a,bertelson.cahen.gutt:1997} 
\[
\Def(M,\omega) \cong \frac{1}{\im \lambda}[\omega] +
\HdR^2(M)[[\lambda]], 
\]
the action $\Phi$ is given by 
\begin{equation}
    \label{eq:mainresult}
    \Phi_{\scriptscriptstyle L}
    ([\omega_{\scriptscriptstyle \lambda}]) 
    = [\omega_{\scriptscriptstyle \lambda}] + 2\pi\im c_1(L),
\end{equation}
where 
$[\omega_{\lambda}]
= (1/{\im \lambda})[\omega] + \sum_{r=0}^\infty[\omega_r] \lambda^r$, 
and $c_1(L)$ is the Chern class of $L$. It immediately follows from
(\ref{eq:mainresult}) that two star products $\star$ and $\star'$
on $M$ are Morita
equivalent if and only if there exists a symplectomorphism
$\psi:M \longrightarrow M$ such that the relative class 
$t(\star, \psi^*(\star'))$ is
$2\pi\im$-integral. The explicit computation of
$\Phi_{\scriptscriptstyle L}$ is based on a local description of
deformed line bundles over $M$, through deformed transition functions,
and the \v{C}ech-cohomological approach to Deligne's relative class
developed in \cite{gutt.rawnsley:1999}. As it turns out, this result
also gives a classification of Hermitian star products on $M$ up to
strong Morita equivalence, a purely algebraic generalization of the
usual notion of strong Morita equivalence of $C^*$-algebras 
\cite{bursztyn.waldmann:2001a,bursztyn.waldmann:2000a:pre}.

By considering star products on cotangent bundles $T^*Q$, we observe
that the integrality condition coming from Morita equivalence can be
interpreted as Dirac's quantization condition for magnetic charges:
We consider the star products $\star_{\kappa}^{-\lambda B}$,
constructed in \cite{bordemann.neumaier.pflaum.waldmann:1998:pre} out
of a $\kappa$-ordered star product $\star_{\kappa}$ on $T^*Q$ and a
magnetic field $B \in \Omega^2(Q)[[\lambda]], dB=0$, and show that
$\star_{\kappa}$ and $\star_{\kappa}^{-\lambda B}$ are Morita
equivalent if and only if $(1/2\pi)B$ is an integral $2$-form.
In this case, well-known $^*$-representations of 
$\star^{-\lambda B}_\kappa$ on sections of  line bundles
\cite{bordemann.neumaier.pflaum.waldmann:1998:pre} are obtained by
means of Rieffel induction of the formal Schr\"{o}dinger
representation of $\star_\kappa$.

After the conclusion of this work, \cite{jurco.schupp.wess:2001:pre}
was brought to our attention;
this paper addresses some related questions and introduces a similar 
local description of quantum line bundles. We note that
(\ref{eq:mainresult}), when written in terms of formal Poisson
structures, coincides with the expression of $\theta'$ in
\cite[pp.~3]{jurco.schupp.wess:2001:pre}. A detailed comparison
between the approaches is in progress.

The paper is organized as follows. In Section~\ref{sec:prelim} we
recall the notions of star products, deformation quantization of
vector bundles and Morita equivalence, and give a local description of
deformed vector bundles in terms of quantum transition matrices,
including Hermitian structures. In Section~\ref{sec:result} we
compute the relative class of Morita equivalent star products on
symplectic manifolds and discuss the main results of the paper. In
Section~\ref{sec:applic} we consider star products on cotangent
bundles and discuss Morita equivalence in terms of Dirac's condition
for magnetic monopoles. We have included two appendices:
Appendix~\ref{sec:explog} recalls some basic facts about
$\star$-exponentials and logarithms used in the paper;
Appendix~\ref{sec:strong} recalls the notions of algebraic Rieffel
induction and strong Morita equivalence. 

\smallskip

\noindent 
{\bf Acknowledgments:} The authors would like to thank Martin Bordemann,
Simone Gutt, Ryszard Nest, Bjorn Poonen and Alan Weinstein for useful
discussions. We also thank Branislav Jur\v{c}o and Peter Schupp for a
valuable discussion clarifying the additional action of the
diffeomorphism group and for bringing
\cite{jurco.schupp.wess:2001:pre} to our attention.

%
%

\section{Preliminaries}
\label{sec:prelim}

\subsection{Star products, deformed vector bundles and Morita equivalence}
\label{sec:star}

Let $(M, \pi_0)$ be a Poisson manifold, where
$\pi_0 \in \Gamma^\infty(\bigwedge^2 TM)$ denotes the Poisson tensor.
The corresponding Poisson bracket is denoted by 
$\{f,g\}:= \pi_0(df,dg)$. 
Let $C^\infty(M)$ be the algebra of complex-valued smooth functions on $M$.
We recall the definition of star products
\cite{bayen.et.al:1978}.
\begin{definition}
    \label{def:starprod}
    A \emph{star product} on a Poisson manifold $(M,\pi_0)$ is a 
    $\mathbb{C}[[\lambda]]$-bilinear associative product on 
    $C^\infty(M)[[\lambda]]$ of the form
    \begin{equation}
        \label{eq:starprod}
        f \star g = \sum_{r=0}^\infty \lambda^r C_r(f,g), 
        \;\; f,g \in C^\infty(M), 
    \end{equation}
    where each $C_r:C^\infty(M)\times C^\infty(M)
    \longrightarrow C^\infty(M)$ is a bidifferential operator,
    $C_0(f, g) = fg$ (pointwise product of
    functions) and $C_1(f, g) - C_1(g, f) = \im \{f, g\}$.
    It is often required that $f \star 1 = f = 1 \star f$. 
\end{definition}
For physical applications, $\lambda$ plays the role of Planck's
constant $\hbar$ as soon as the convergence of \eqref{eq:starprod} can
be established. The existence of star products on symplectic manifolds
was proven in
\cite{dewilde.lecomte:1983b,fedosov:1994a,omori.maeda.yoshioka:1991};
for arbitrary Poisson manifolds, it follows from Kontsevich's
formality theorem \cite{kontsevich:1997:pre}.

Two star products $\star$ and $\star'$ are called \emph{equivalent} if
there exist differential operators 
$T_r:C^\infty(M) \longrightarrow C^\infty(M)$ so that  
$T = \id + \sum_{r=1}^\infty \lambda^r T_r$ satisfies
\begin{equation}
    \label{eq:equiv}
    T(f\star'g) = T(f) \star T(g), \;\; f,g \in C^\infty(M).
\end{equation}
The equivalence class of a star product $\star$ on $(M,\pi_0)$ will be
denoted by $[\star]$. We let
\begin{equation}
    \label{eq:moduli}
    \Def(M,\pi_0)
    :=
    \{ [\star], \; \star\; \mbox{ a star product on } (M,\pi_0)\}.
\end{equation}
For symplectic manifolds, the moduli space (\ref{eq:moduli}) admits a
cohomological description
\cite{nest.tsygan:1995a,bertelson.cahen.gutt:1997} that will be
recalled in Section~\ref{sec:relative}.
We note that the group of Poisson diffeomorphisms of $M$ acts
naturally on $\Def(M,\pi_0)$: $\star' = \psi^*(\star)$ is defined
by $\psi^*(f\star'g) = \psi^*f \star \psi^*g$.

A classical result of Serre and Swan \cite[Chap.~XIV]{bass:1968} 
asserts that finite dimensional complex vector bundles over $M$
naturally correspond to finitely generated projective modules over
$C^\infty(M)$ (with equivalence functor 
$E \mapsto \Gamma^\infty(E)$). This motivates the following definition
\cite[Def.~3.1]{bursztyn.waldmann:2000b}: Let $E \to M$ be a
$k$-dimensional complex vector bundle, and let $\star$ be a star
product on $M$. 
\begin{definition}
    \label{def:defvecbund}
    A \emph{deformation quantization} of $E \to M$ with respect to
    $\star$ is a $\mathbb{C}[[\lambda]]$-bilinear map
    $\bullet: \Gamma^\infty(E)[[\lambda]] 
    \times C^\infty(M)[[\lambda]] \longrightarrow
    \Gamma^\infty(E)[[\lambda]]$ 
    satisfying
    $s \bullet (f \star g) = (s \bullet f) \bullet g$ and
    so that
    \begin{equation}
        \label{eq:bullet}
        s \bullet f = \sum_{r=0}^\infty \lambda^r R_r (s, f),
    \end{equation}
    where each 
    $R_r: \Gamma^\infty(E)\times C^\infty(M) 
    \longrightarrow \Gamma^\infty(E)$
    is bidifferential and $R_0(s,f)=sf$ (pointwise multiplication of
    sections by functions).
\end{definition}
Two deformations $\bullet$ and $\bullet'$ are called \emph{equivalent}
if there exist differential operators 
$T_r: \Gamma^\infty(E) \longrightarrow \Gamma^\infty(E)$ so that
$T = \id + \sum_{r=1}^\infty \lambda^r T_r$ satisfies
\begin{equation}
    \label{eq:equivmod}
    T(s \bullet' f) = (Ts) \bullet f, \;\; s \in \Gamma^\infty(E), \; 
    f \in C^\infty(M).
\end{equation}
The following result was proven in
\cite[Prop.~2.6]{bursztyn.waldmann:2000b}.
\begin{proposition}
    \label{prop:existence}
    Let $E \to M$ be a vector bundle, and let $\star$ be a star
    product on $M$. Then there exists a deformation quantization
    $\bullet$ of $E$ with respect to $\star$, which is unique up to
    equivalence. The right module 
    $(\Gamma^\infty(E)[[\lambda]], \bullet)$ is finitely generated and
    projective over $(C^\infty(M)[[\lambda]], \star)$, and any
    finitely generated projective module over this algebra arises in
    this way.
\end{proposition}

Let $\mathcal{E} = \Gamma^\infty(E)$, considered as a right
$C^\infty(M)$-module, and 
$\boldsymbol{\mathcal{E}} = (\mathcal{E}[[\lambda]], \bullet)$,
considered as a right $(C^\infty(M)[[\lambda]],\star)$-module. We
recall that $\End(\mathcal{E}) \cong \Gamma^\infty(\End E)$, and
$\End(\boldsymbol{\mathcal{E}})$ is isomorphic to
$\Gamma^\infty(\End E)[[\lambda]]$ as a $\mathbb{C}[[\lambda]]$-module
\cite[Cor.~2.4]{bursztyn.waldmann:2000b}.

If $E=L \to M$ is a complex line bundle, 
then $\Gamma^\infty(\End(L))\cong C^\infty(M)$, and any 
$\mathbb{C}[[\lambda]]$-module isomorphism 
$\Gamma^\infty(\End(L))[[\lambda]] \longrightarrow
\End(\boldsymbol{\mathcal{E}})$ 
determines a new star product $\star'$ on $M$. 
We choose this isomorphism  so that the corresponding left action
of $(C^\infty(M)[[\lambda]],\star')$ on $\boldsymbol{\mathcal{E}}$
deforms the multiplication of sections by functions.
As shown in 
\cite[Lem.~3.4]{bursztyn:2001:pre}, $\star$ and $\star'$ correspond to
the same Poisson structure on $M$, and this procedure gives rise to a
well-defined map 
\begin{equation}
    \label{eq:bijec}
    \Phi_{\scriptscriptstyle L}: 
    \Def(M,\pi_0) \longrightarrow \Def(M,\pi_0), 
    \;\;\; [\star] \mapsto [\star'].
\end{equation}
It is simple to check that the map $\Phi_{\scriptscriptstyle L}$
depends only on the isomorphism class of $L$, denoted by 
$[L] \in \Pic(M) \cong \CH^2(M, \mathbb{Z})$, where $\Pic(M)$ is the
Picard group of $M$. The following result was proven in
\cite[Thm.~4.1]{bursztyn:2001:pre}.
\begin{proposition}
    \label{prop:action}
    The map 
    $\Phi:\Pic(M)\times \Def(M,\pi_0)  
    \longrightarrow \Def(M,\pi_0)$,
    $([L],[\star]) \mapsto 
    \Phi_{\scriptscriptstyle L}([\star])$ 
    defines an action of $\Pic(M)$ on the set of equivalence classes
    of star products on $M$, and two star products
    $\star$ and $\star'$ on $M$ are Morita
    equivalent if and only if there exists a Poisson diffeomorphism
    $\psi: M \longrightarrow M$ such that
    the classes $[\star]$ and $[\psi^*(\star)]$ lie in the same $\Phi$-orbit.
\end{proposition}

Recall that two unital algebras $\mathcal A, \mathcal{B}$ (over some
ground ring $R$) are called \emph{Morita equivalent} if they have
equivalent categories of left modules \cite{jacobson:1989};
alternatively, there must exist a full finitely generated projective
right $\mathcal A$-module 
$\mathcal{E}_{\scriptscriptstyle\mathcal{A}}$ so that 
$\mathcal B \cong \End(\mathcal{E}_{\scriptscriptstyle\mathcal{A}})$. 
The bimodule 
${}_{\scriptscriptstyle{\mathcal{B}}}\mathcal{E}_{\scriptscriptstyle{\mathcal{A}}}$
is called a  \emph{$(\mathcal B, \mathcal A)$-equivalence bimodule}.

The Picard group of a unital $R$-algebra $\mathcal A$,
$\Pic(\mathcal A)$, is defined as the set of isomorphism classes of
$(\mathcal A,\mathcal A)$-equivalence bimodules, with group operation
given by tensor product. If $\mathcal A = C^\infty(M)$, then the
algebraic Picard group $\Pic(C^\infty(M))$ can be identified with the
geometric Picard group $\Pic(M)$.  

Let $\star$ be a star product on $M$, and 
$\boldsymbol{\mathcal A}=(C^\infty(M)[[\lambda]],\star)$. We note that the
isotropy  group of $\Phi$ at $[\star]$ can be identified with a subgroup of
$\Pic(\boldsymbol{\mathcal A})$. In Section~\ref{sec:relative}, 
we will give an explicit description of $\Pic(\boldsymbol{\mathcal A})$
for certain star-product algebras on symplectic manifolds.
We also describe the orbit space $\Def(M,\pi_0)/{\Pic(M)}$
for $\pi_0$ symplectic.


\subsection{A local description of deformed vector bundles}
\label{sec:local}

Let $E \to M$ be a $k$-dimensional smooth complex vector bundle over a
smooth manifold $M$, and let 
$\{\mathcal{O}_\alpha \}$ be a good cover of $M$. 
Let us fix $\{e_{\alpha,i}\}, i=1 \ldots k$, basis of
$\Gamma^\infty(E |_{\mathcal{O}_\alpha})$, and let 
$e_\alpha = (e_{\alpha,1},\ldots,e_{\alpha,k})$ be the corresponding
frame. Such a choice defines trivialization maps 
$\psi_\alpha: \Gamma^\infty(E |_{\mathcal{O}_\alpha}) 
\longrightarrow C^\infty(\mathcal{O}_\alpha)^k$.
On overlaps $\mathcal{O}_\alpha \cap \mathcal{O}_\beta$, we define
transition matrices 
$\phi_{\alpha \beta} = \psi_{\alpha}\psi_{\beta}^{-1} 
\in M_k(C^\infty(\mathcal{O}_\alpha \cap \mathcal{O}_\beta))$. 
Clearly $\phi_{\alpha \beta} = \phi_{\beta \alpha}^{-1}$, and on
triple intersections we have 
\begin{equation}
    \label{eq:cocycle}
    \phi_{\alpha \beta}\phi_{ \beta \gamma}\phi_{\gamma \alpha} = 1.
\end{equation}
We will see that similar constructions can be carried out for deformed
vector bundles (see also \cite{jurco.schupp.wess:2000:pre,jurco.schupp.wess:2001:pre}).
Let $\star$ be a star product on $M$, and let us fix a deformation
$\boldsymbol{\mathcal E}=(\mathcal{E}[[\lambda]],\bullet)$,
$\mathcal E = \Gamma^\infty(E)$, with respect to $\star$. A simple
induction shows the following result.
\begin{lemma}
    \label{lemma:DefBaseSections}
    Let
    $\qe_{\alpha,i}
    = e_{\alpha,i} + \lambda e^{(1)}_{\alpha,i} + \cdots \in 
    \Gamma^\infty(E|_{\mathcal{O}_\alpha})[[\lambda]]$ be arbitrary
    deformations of the classical bases sections $e_{\alpha,i}$. Then
    for any global section $\qs \in \Gamma^\infty(E)[[\lambda]]$ there
    exist unique local functions
    $\qs^i_\alpha \in C^\infty(\mathcal{O}_\alpha)[[\lambda]]$
    such that
    \begin{equation}
        \label{eq:LocalDefBase}
        \qs|_{\mathcal{O}_\alpha} 
        = \sum_{i=1}^k\qe_{\alpha, i} \bullet \qs^i_\alpha.
    \end{equation}
\end{lemma}
We shall write $\qe_\alpha = (\qe_{\alpha,1}, \ldots, \qe_{\alpha,k})$
for the deformed frame, and 
$\qs_\alpha = (\qs_{\alpha}^1, \ldots, \qs_\alpha^k)$ for the deformed
coefficient functions of a section $s$. As in the case of ordinary
vector bundles, (\ref{eq:LocalDefBase}) induces 
$\mathbb{C}[[\lambda]]$-linear trivialization isomorphisms
$\qpsi_\alpha: \Gamma^\infty(E|_{\mathcal{O}_\alpha})[[\lambda]]
\to C^\infty(\mathcal{O}_\alpha)^k[[\lambda]]$,
\begin{equation}
    \label{eq:DefTrivComponents}
    \qpsi_\alpha=(\qpsi_\alpha^1,\ldots,\qpsi_\alpha^k),\;\;\;
    \qpsi_\alpha^i (\qs)
    = \qpsi_\alpha^i (\sum_j \qe_{\alpha,j} \bullet \qs_\alpha^j)
    = \qs_\alpha^i, 
\end{equation}
satisfying 
\begin{equation}
    \label{eq:DefTrivModMorph}
    \qpsi_\alpha(\qs \bullet f) = \qpsi_\alpha (\qs) \star f, \;\;
    \mbox{ for } f \in C^\infty(M).
\end{equation}
Clearly $\qpsi_\alpha (\qs) = \qs_\alpha$. It is simple to check that
$\qpsi_\alpha$ deforms $\psi_\alpha$, i.e., 
$\qpsi_\alpha= \psi_\alpha \mbox{ mod } \lambda$.

On overlaps $\mathcal{O}_\alpha \cap \mathcal{O}_\beta$, we define
deformed transition matrices 
\begin{equation}
    \label{eq:DefTransMatrix}
    \qphi_{\alpha \beta} = 
    \qpsi_{\alpha}\circ \qpsi_{\beta}^{-1} \in
    M_k(C^\infty(\mathcal{O}_\alpha\cap\mathcal{O}_\beta))[[\lambda]],   
\end{equation}
satisfying $\qs_\alpha = \qphi_{\alpha\beta} \star \qs_\beta$.
We note that $\qphi_{\alpha \beta}=\qphi_{\beta \alpha} ^{-1}$ (with
respect to $\star$), and the following deformed cocycle
condition holds:
\begin{equation}
    \label{eq:QCocycleMatrix}
    \qphi_{\alpha\beta} \star
    \qphi_{\beta\gamma} \star
    \qphi_{\gamma\alpha} = \Unit.
\end{equation}

If $\qA \in \End(\boldsymbol{\mathcal E})$, then it is locally
represented by a matrix 
$\qA_\alpha \in M_k(C^\infty(\mathcal{O}_\alpha))[[\lambda]]$ 
satisfying $\qA(\qs)_\alpha = \qA_\alpha \star \qs_\alpha$. 
On overlaps $\mathcal{O}_\alpha \cap \mathcal{O}_\beta$, we have
\begin{equation}
    \label{eq:ChangeQuantACoeff}
    {\qA_\beta}
    = {\qphi_{\beta\alpha}}
    \star {\qA_\alpha} \star
    {\qphi_{\alpha\beta}}.
\end{equation}
As in the classical case, a collection
$\{\qA_\alpha\}$, 
$\qA_\alpha \in M_k(C^\infty(\mathcal{O}_\alpha))[[\lambda]]$,
satisfying (\ref{eq:ChangeQuantACoeff}) determines a global
endomorphism of the deformed bundle. It is simple to see that the
composition of endomorphisms corresponds locally to the deformed
product of matrices:
\begin{equation}
    \label{eq:EndoProd}
    (\qA\circ \qB)_{\alpha} = \qA_\alpha \star \qB_\alpha.
\end{equation}
\begin{remark}
    \label{remark:localequiv}
    One can define an explicit $\mathbb{C}[[\lambda]]$-module
    isomorphism 
    $T: \Gamma^\infty(\End(E))[[\lambda]] \to
    \End(\boldsymbol{\mathcal{E}})$ 
    by patching local maps as follows. Let $\{\chi_\alpha\}$ be a
    quadratic partition of unity subordinated to
    $\{\mathcal{O}_\alpha\}$ (i.e. 
    $\supp\chi_\alpha \subseteq \mathcal{O}_\alpha$, and
    $\sum_\alpha \cc{\chi_\alpha} \chi_\alpha = 1$). Then 
    \begin{equation}
        \label{eq:localexpress}
        T_\alpha(A) = \sum_\gamma
        {\qphi_{\alpha\gamma}}
        \star \cc{\chi_\gamma} \star A_\gamma \star \chi_\gamma \star
        {\qphi_{\gamma\alpha}}
    \end{equation}
    is well defined on $\mathcal{O}_\alpha$. Here $A_\alpha$ are the
    local matrices of $A \in \Gamma^\infty(\End(E))[[\lambda]]$ with
    respect to the undeformed trivialization maps
    $\psi_{\alpha}$. The collection $\{T_\alpha\}$ satisfies
    condition (\ref{eq:ChangeQuantACoeff}), and hence defines the
    desired global map $T$. In lowest order $T_\alpha(A)$
    just reproduces $A_\alpha$.
\end{remark}

\subsection{Hermitian structures}
\label{sec:HermBundle}

For completeness, we will briefly indicate how deformed
Hermitian structures \cite{bursztyn.waldmann:2000b} can be treated
locally. In this section, $\star$ will be a Hermitian star product on
$M$, i.e. $\cc{f\star g} = \cc{g}\star \cc{f}$.

Let $E \to M$ be equipped with a Hermitian fiber metric
$h_0$. A \emph{deformation quantization} of $h_0$ with respect to a
deformation $\bullet$ of $E$ is a
$C^\infty(M)[[\lambda]]$-valued Hermitian inner product $\qh$ on
$\Gamma^\infty(E)[[\lambda]]$ 
(see Definition \ref{def:prehilbert}) such that  
\begin{equation}
    \label{eq:HermDef}
    \qh(\qs, \qs') = \sum_{r=0}^\infty \lambda^r h_r (\qs, \qs')
\end{equation}
with bidifferential operators $h_r:\Gamma^\infty(E)\times\Gamma^\infty(E)
\longrightarrow C^\infty(M)$.

Let $\boldsymbol{\mathcal E}$ denote the
$(C^\infty(M)[[\lambda]],\star)$-module 
$(\Gamma^\infty(E)[[\lambda]],\bullet)$.
Two deformations $\qh$ and $\qh'$ are called \emph{isometric} if there
exists a module isomorphism 
\begin{equation}
    \label{eq:isometry}
   \boldsymbol{U} = \id + \sum_{r=1}^\infty \lambda^r U_r:
    \boldsymbol{\mathcal E} \longrightarrow \boldsymbol{\mathcal E},
\end{equation}
with differential
operators $U_r:\Gamma^\infty(E)\longrightarrow \Gamma^\infty(E)$, so that
\begin{equation}
    \label{eq:isometric}
    \qh( \boldsymbol{U}\qs, \boldsymbol{U}\qs') = \qh' (\qs, \qs')
\end{equation}
for all $\qs, \qs' \in \Gamma^\infty(E)[[\lambda]]$. From
\cite{bursztyn.waldmann:2000b} we have the following result.
\begin{lemma}
    \label{lemma:HermDefExist}
    Let $E \to M$ be a vector bundle with Hermitian fiber metric $h_0$,
    and let $\bullet$ be a deformation quantization of $E$. Then there
    exists a deformation quantization $\qh$ of $h_0$ and any two such
    deformations are isometric.
\end{lemma}

Let $\qh$ be a deformation of $h_0$. We can construct local
orthonormal frames $\qe_\alpha$ with respect to $\qh$:
\begin{lemma}
    \label{lemma:OrthoFrames}
    Let $\widetilde{\qe}_\alpha$ be a local frame for
    $\Gamma^\infty(E)[[\lambda]]$ such that the zeroth order is an
    orthonormal frame with respect to $h_0$. Then there exists a
    matrix 
    $\qV = \id + \sum_{r=1}^\infty \lambda^r V_r \in
    M_k(C^\infty(\mathcal{O}_\alpha)[[\lambda]])$ 
    such that $\qe_\alpha := \widetilde{\qe}_\alpha \bullet \qV$ is an
    orthonormal frame with respect to $\qh$, i.e. one has
    \begin{equation}
        \label{eq:OrthoFrame}
        \qh(\qe_{\alpha,i}, \qe_{\alpha,j}) = \delta_{ij}.
    \end{equation}
\end{lemma}
\begin{proof}
    Let $\qH$ be the Hermitian matrix defined by 
    $\qH_{ij} =
    \qh(\widetilde{\qe}_{\alpha,i},\widetilde{\qe}_{\alpha,j})$.
    Then $\qH = \id + \sum_{r=1}^\infty \lambda^r H_r$, since the
    zeroth order of $\widetilde{\qe}_\alpha$ is orthonormal with
    respect to $h_0$. From \cite[Lem.~2.1]{bursztyn.waldmann:2000b} we
    know that there exists a matrix 
    $\qU = \id +  \sum_{r=1}^\infty \lambda^r U_r$ such that
    $\qU^* \star \qU = \qH$. Then $\qV = \qU^{-1}$ is the desired
    transformation. 
   
\end{proof}

Hence we can always assume that we have local orthonormal
frames $\qe_\alpha$ on each patch $\mathcal{O}_\alpha$. Obviously,
the transition functions are unitary in this case:
\begin{lemma}
    \label{lemma:UnitaryTrans}
    Let $\{\qe_\alpha\}_{\alpha \in I}$ be local orthonormal
    frames. Then we have 
    \begin{equation}
        \label{eq:UnitrayTrans}
        \qphi_{\alpha\beta}^* 
        = \qphi_{\beta\alpha} 
        = \qphi_{\alpha\beta}^{-1}
    \end{equation}
    and $\qh(\qs, \qs') = \SP{\qs_\alpha, \qs'_\alpha}$ is just the
    canonical Hermitian inner product on 
    $C^\infty(\mathcal{O_\alpha})[[\lambda]]^k$ for the coefficient
    functions . If $\qA \in \End(\boldsymbol{\mathcal E})$, then the
    local matrices of $\qA$ and $\qA^*$ are related by
    \begin{equation}
        \label{eq:AAstarLocal}
        (\qA_\alpha)^* = (\qA^*)_\alpha.
    \end{equation}
\end{lemma}
Note that, in this case, the isomorphism (\ref{eq:localexpress}) is
compatible with the $^*$-structures.

%
%

\section{Morita equivalent star products on symplectic manifolds}
\label{sec:result}

\subsection{Deligne's relative class (after Gutt and Rawnsley)}
\label{sec:deligne}

Let $(M,\omega)$ be a symplectic manifold. In this case, it was shown
in \cite{nest.tsygan:1995a,bertelson.cahen.gutt:1997} that there exists a
bijection  
\begin{equation}
    \label{eq:charac}
    c: \Def(M,\omega) \longrightarrow 
    \frac{1}{\im \lambda}[\omega] + \HdR^2(M)[[\lambda]], 
\end{equation}
characterizing the moduli space of equivalence classes of star
products on $M$ in cohomological terms. For a star product $\star$,
$c(\star)$ is called its \emph{characteristic class}. A
\v{C}ech-cohomological description of these characteristic classes can
be found in \cite{gutt.rawnsley:1999,deligne:1995}.

For two star products $\star$, $\star'$ on $M$, their \emph{relative
  class} is defined by
\begin{equation}
    t(\star',\star) = c(\star')-c(\star) \in \HdR^2(M)[[\lambda]].
\end{equation}
A purely \v{C}ech-cohomological construction of $t(\star',\star)$ was
given in \cite{gutt.rawnsley:1999}, and we will briefly recall it.

Let us fix a good cover $\{\mathcal{O}_\alpha\}$ of $M$ and
star products $\star, \star'$. Then any two star products are
equivalent on $\mathcal{O}_\alpha$, see
e.g.~\cite[Cor.~3.2]{gutt.rawnsley:1999}. Thus, for each $\alpha$, we
can find an equivalence transformation between $\star$ and $\star'$,
$T_\alpha = \id + \sum_{r=1}^\infty \lambda^r T_\alpha^{(r)}$, where
each $T_\alpha^{(r)}$ is a differential operator on
$C^\infty(\mathcal{O}_\alpha)$. On the overlap
$\mathcal{O}_\alpha \cap \mathcal{O}_\beta$,
the map $T_\alpha^{-1} \circ T_\beta$ is a $\star$-automorphism
starting with the identity. Since 
$\mathcal{O}_\alpha \cap \mathcal{O}_\beta$ is contractible,
the automorphism $T_\alpha^{-1} \circ T_\beta$ is inner, and therefore
there exists a function 
$t_{\alpha\beta} \in 
C^\infty (\mathcal{O}_\alpha \cap \mathcal{O}_\beta)[[\lambda]]$
such that (see Prop.~\ref{proposition:Exp})
\begin{equation}
    \label{eq:TalphaTbeta}
    T_\alpha^{-1} \circ T_\beta (f)
    = \eu^{[t_{\alpha\beta}, \cdot]}(f)
    = \Exp(t_{\alpha\beta}) \star f \star \Exp(-t_{\alpha\beta}).
\end{equation}
Since
\begin{equation}
    \label{eq:ProdTabc}
    T_\alpha^{-1} \circ T_\beta
    \circ T_\beta^{-1} \circ T_\gamma \circ
    T_\gamma^{-1} \circ T_\alpha = \id,
\end{equation}
the element 
$\Exp(t_{\alpha\beta}) \star \Exp(t_{\beta\gamma}) \star
\Exp(t_{\gamma\alpha})$ 
must be central. Thus 
\begin{equation}
    \label{eq:TabsDefAgain}
    t_{\alpha\beta\gamma} =
    t_{\alpha\beta} \circ_\star t_{\beta\gamma} \circ_\star
    t_{\gamma\alpha}
    \in \mathbb{C}[[\lambda]]
\end{equation}
defines a \v{C}ech cochain on $M$ with values in
$\mathbb{C}[[\lambda]]$. This cochain turns out to be a cocycle
\cite{gutt.rawnsley:1999}, and the \v{C}ech class
$[t_{\alpha\beta\gamma}]$ (viewed as a class in
$\HdR^2(M)[[\lambda]]$) is the relative class $t(\star',\star)$.

\subsection{The relative class of Morita equivalent star products}
\label{sec:relative}

We will now use the results in Sections~\ref{sec:deligne} and
\ref{sec:local} to compute the relative class of two Morita equivalent
star products on a symplectic manifold $(M,\omega)$, providing an
explicit description of the orbit space $\Def(M,\omega)/\Pic(M)$.
\begin{theorem}
    \label{thm:main}
    Let $L \to M$ be a complex line bundle over a symplectic manifold
    $M$. Suppose $\star$, $\star'$ are star products on $M$, with
    $\Phi_{\scriptscriptstyle L}([\star])=[\star']$. 
    Then $t(\star',\star)= 2\pi\im c_1(L)$,
    where $c_1(L)$ is the Chern class of $L$.
\end{theorem}
\begin{proof}
Let $\{\mathcal{O}_\alpha\}$ be a good cover of $M$, and let us fix
deformed trivialization maps $\qpsi_\alpha$ and transition functions
$\qphi_{\alpha\beta}$ as in Section~\ref{sec:local}.
Let $\boldsymbol{\mathcal E}=(\Gamma^\infty(L)[[\lambda]],\bullet)$ be
a deformation quantization of $L$ with respect to $\star$. Let 
$T:(C^\infty(M)[[\lambda]],\star') 
\longrightarrow \End(\boldsymbol{\mathcal E})$ 
be a $\mathbb{C}[[\lambda]]$-algebra isomorphism, that, by
\cite{bursztyn.waldmann:2000b}, can be chosen to preserve
supports (see Remark~\ref{remark:localequiv}).
Such a $T$ gives rise to a collection of local maps
\[
T_\alpha: C^\infty(\mathcal{O}_\alpha)[[\lambda]]
\longrightarrow  C^\infty(\mathcal{O}_\alpha)[[\lambda]],
\]
by $T_\alpha(f) = T(f)_\alpha$, satisfying 
$T_\alpha = \id \mbox{ mod }\lambda$
and $T_\alpha f \star T_\alpha g = T_\alpha(f\star'g)$ 
(by (\ref{eq:EndoProd})). It follows from (\ref{eq:ChangeQuantACoeff})
that
$T_\beta(f) = 
\qphi_{\beta\alpha} \star T_\alpha(f) \star \qphi_{\alpha\beta}$, 
and therefore
\begin{equation}
    \label{eq:localequivexp}
    T_\alpha^{-1}T_\beta(f)
    = \qphi_{\alpha\beta}\star f \star\qphi_{\beta\alpha}.
\end{equation}
Since $\phi_{\alpha \beta}$ is invertible and 
$\qphi_{\alpha \beta} = \phi_{\alpha \beta}\mbox{ mod }\lambda$, we
can write (see Appendix \ref{sec:explog})
\[
\qphi_{\alpha \beta} = \Exp(t_{\alpha \beta}),
\]
for some 
$t_{\alpha\beta} = t_{\alpha\beta}^{(0)}
+ \sum_{r=1}^\infty \lambda^r t_{\alpha\beta}^{(r)}
\in C^\infty(\mathcal{O}_\alpha \cap \mathcal{O}_\beta)[[\lambda]]$, 
and $\phi_{\alpha \beta}= \eu^{ t_{\alpha\beta}^{(0)}}$.

The deformed cocycle condition (\ref{eq:QCocycleMatrix}) and
Prop.~\ref{proposition:Exp} imply that, on triple intersections
$\mathcal{O}_\alpha\cap \mathcal{O}_\beta \cap \mathcal{O}_\beta$, the
function 
$t_{\alpha\beta\gamma} := t_{\alpha\beta}
\circ_\star t_{\beta\gamma} \circ_\star t_{\gamma\alpha}$ 
must satisfy
\[
t_{\alpha\beta\gamma} = 2 \pi \im n_{\alpha\beta\gamma}, \,
\; \textrm{ with }
n_{\alpha\beta\gamma} \in \mathbb{Z}.
\]    
This shows that $\frac{1}{2\pi\im}t(\star',\star)$ is integral and
does not depend on $\lambda$. Since the classical limit of
$\circ_\star$ is just the usual addition, we get 
\[
t_{\alpha \beta \gamma}= t_{\alpha\beta}^{(0)} + t_{\beta\gamma}^{(0)}
+ t_{\gamma\alpha}^{(0)} = 2 \pi \im n_{\alpha\beta\gamma}.
\]
But the complex \v{C}ech class defined by 
$\frac{1}{2\pi\im}(t_{\alpha\beta}^{(0)} + t_{\beta\gamma}^{(0)} +
t_{\gamma\alpha}^{(0)})$, viewed as a de Rham class, is the Chern class of
$L$. Thus $t(\star',\star)=2\pi\im c_1(L)$.

\end{proof}

Let $\HdR^2(M,\mathbb Z)$ denote the image of the usual
map $i:\CH^2(M, \mathbb{Z}) \longrightarrow \HdR^2(M,\mathbb C)$.
\begin{corollary}
    \label{corollary:StarME}
    Two star products $\star, \star'$ on a symplectic manifold $M$ are
    Morita equivalent if and only if there exists a symplectomorphism
    $\psi: M \longrightarrow M$ such that
    \begin{equation}
        \label{eq:StarME}
        \frac{1}{2\pi \im}t(\star',\psi^*(\star))\in \HdR^2(M,\mathbb Z).
    \end{equation}
\end{corollary}
An immediate consequence of Theorem \ref{thm:main} is the following
explicit expression for the action $\Phi$ in terms of the
characteristic classes of star products:
\begin{equation}
    \label{eq:explicitPhi}
    \Phi_{\scriptscriptstyle L}([\omega_\lambda])
    = [\omega_\lambda]+ 2\pi \im c_1(L),
\end{equation}
where 
$[\omega_\lambda] = 
(1/{\im \lambda})[\omega]+ \sum_{r=0}^\infty[\omega_r]\lambda^r$.
The orbit space $\Def(M,\omega)/\Pic(M)$ is just a trivial fibration
over the torus $\HdR^2(M,\mathbb C)/\HdR(M,\mathbb Z)$, with fiber 
$\HdR^2(M,\mathbb C)[[\lambda]]$.

It is clear that the isotropy group of $\Phi$, for any 
$[\star] \in \Def(M,\omega)$, is isomorphic to the subgroup
$\mathcal{T}(M) :=\{[L]\in \Pic(M), \; c_1(L)=0\}\subseteq \Pic(M)$ of
\emph{flat line bundles}. 

Let  $\star$ be  a star product on $(M,\omega)$
with $c(\star) = [\omega]/{\im \lambda} + O(\lambda)$ (i.e.,
$c_0(\star)=0$). Since $c(\psi^*(\star)) = \psi^*c(\star)$, it follows that
in this case $\Pic(C^\infty(M)[[\lambda]],\star)$ is isomorphic to
the isotropy group of $\Phi$ at $[\star]$. Hence we have

\begin{corollary}
    \label{corollary:defPic}
    Let $\star$ be a star product on $(M,\omega)$ with
    $c_0(\star)=0$. Then the Picard
    group of the algebra $(C^\infty(M)[[\lambda]],\star)$ is
    isomorphic to $\mathcal{T}(M)$.
\end{corollary}
Under the usual identification $\Pic(M)\cong \CH^2(M, \mathbb{Z})$,
$\mathcal{T}(M)$ correspond to torsion elements in 
$\CH^2(M, \mathbb{Z})$. Hence if $\CH^2(M, \mathbb{Z})$ is free,
$\Phi$ is faithful and the Picard groups of the deformed algebras 
with $c_0(\star)=0$ are
trivial.
\begin{corollary}
    \label{corollary:flatlines}
    Let $L \to M$ be a line bundle over $(M, \omega)$. Then
    $\Gamma^\infty(L)[[\lambda]]$ has a $\star$-bimodule structure
    deforming the classical one if and only if $L$ is flat.
\end{corollary}

\subsection{Strong Morita equivalence of star products}
\label{sec:sme}

We now observe that Theorem~\ref{thm:main} also provides a complete
classification of Hermitian star products up to strong Morita
equivalence, see Appendix~\ref{sec:strong}. The following lemma
should be well-known.
\begin{lemma}
    \label{lemma:ExpDAuto}
    Let $\mathcal{A}$ be a $k$-algebra, where $k$ is a commutative ring
    with $\mathbb{Q} \subseteq k$. Let 
    $D$ and $T = \exp(\lambda D)$ be $k[[\lambda]]$-module
    endomorphisms of $\mathcal{A}[[\lambda]]$.
    If $\star$ is a formal
    associative deformation for $\mathcal{A}$, then
    $T$ is a $\star$-automorphism if and only if $D$ is a $\star$-derivation.    
\end{lemma}
\begin{proof}
    If $D$ is a $\star$-derivation, then  $T$ is clearly a
    $\star$-automorphism. For the converse, define  
    $E(a,b) = D(a \star b) - Da \star b - a \star Db$. It follows that
    \[
    D^k(a\star b) = \sum_{l=0}^k \binom{k}{l} D^la \star D^{k-l}b
    + \sum_{r, s, t=0}^{k-1} c_{rst}^{(k)} \, D^r E(D^sa, D^tb)
    \]
    with some rational coefficients $c^{(k)}_{rst}$,
    obtained by recursion. From the fact that $T$ is an automorphism,
    we obtain
    \[
    E(a,b) = - \lambda \sum_{k=2}^\infty \frac{\lambda^{k-2}}{k!} 
    \sum_{r, s, t =0}^{k-1} c_{rst}^{(k)} \, D^r E(D^sa, D^tb).
    \]
    This equation can be seen as a fixed point condition for a
    $k[[\lambda]]$-linear operator acting on $k[[\lambda]]$-bilinear maps
    on $\mathcal{A}[[\lambda]]$, and this operator is
    clearly contracting in the $\lambda$-adic topology. Thus, by Banach's fixed 
    point theorem,
    there exists a unique
    fixed point, which must be
    $0$ (see
    e.g.~\cite[App.~A]{bordemann.neumaier.waldmann:1999}). Therefore $E = 0$,
    and $D$ is a derivation. 
    
\end{proof}
\begin{corollary}
    \label{corollary:StarEquiv}
    Let $\star$, $\star'$ be Hermitian star products on a Poisson manifold $M$.
    Then $\star$ is equivalent to $\star'$ if and
    only if $\star$ is $^*$-equivalent to $\star'$.    
\end{corollary}
\begin{proof}
    Let $T = \id + \sum_{r=1}^\infty \lambda^r T_r$ be an equivalence,
    $T(f \star g) = Tf \star' Tg$. Then $f^\dag := T^{-1}(\cc{Tf})$
    defines a new $^*$-involution for $\star$ of the form
    $f^\dag = S \cc f$, where 
    $S = \id + \sum_{r=1}^\infty \lambda^r S_r$ is a $\star$-automorphism. 
    We can write
    $S = \eu^{\im\lambda D}$, where $D$ is a real derivation of
    $\star$. 
    Thus $S^{1/2}$ is still a $\star$-automorphism, and the map 
    $U = T S^{1/2}$ is a $^*$-equivalence between $\star$ and $\star'$.
    
\end{proof}
\begin{theorem}
    \label{theorem:StrongME}
    Let $\star$ and $\star'$ be Hermitian star products on a Poisson manifold
    $M$. Then $\star$ and $\star'$ are strongly Morita equivalent if
    and only if they are Morita equivalent.
    \end{theorem}
\begin{proof}
    Assume that $\star$ and
    $\star'$ are Morita equivalent via a line bundle $L$. Equip
    $L$ with a Hermitian fiber metric $h_0$, and let 
    $\boldsymbol{\mathcal E}
    = (\Gamma^\infty(L)[[\lambda]],\bullet, \boldsymbol{h})$ 
    be a quantization respect to $\star$. The endomorphisms
    $\End(\boldsymbol{\mathcal E})$ form a $^*$-algebra strongly
    Morita equivalent to $(C^\infty(M)[[\lambda]], \star)$, see
    \cite{bursztyn.waldmann:2000b}. This algebra is isomorphic to
    $(C^\infty(M)[[\lambda]], \star')$, and, by
    Lemma~\ref{corollary:StarEquiv}, we  
    can chose the isomorphism to be a $^*$-isomorphism. Hence $\star$
    and $\star'$ are strongly Morita equivalent. For the converse,
    see \cite[Sec.~7]{bursztyn.waldmann:2001a}.
    
\end{proof}
\begin{corollary}\label{cor:strongME} 
    If $M$ is symplectic, and $\star, \star'$ are Hermitian star
    products, then they are strongly Morita equivalent if and only if
    there exists a symplectomorphism $\psi:M \longrightarrow M$ such that
    $c(\psi^*(\star'))-c(\star)$ is $2\pi\im$-integral.
\end{corollary}
We note that a similar result holds for $C^*$-algebras
\cite{beer:1982}: two unital $C^*$-algebras are strongly Morita
equivalent if and only if they are Morita equivalent as unital rings.

%
%

\section{Application}
\label{sec:applic}
In this section we shall consider star products on cotangent bundles
$\pi:T^*Q\to Q$, motivated by the importance of this class of
symplectic manifolds in physical applications.

\subsection{Star products on $T^*Q$}
\label{sec:StarTQ}

We will briefly recall the construction of star products on 
cotangent bundles in order to set up our notation. The reader is
referred to  
\cite{bordemann.neumaier.waldmann:1998,bordemann.neumaier.waldmann:1999,bordemann.neumaier.pflaum.waldmann:1998:pre}
for details.

For $\gamma \in \Gamma^\infty (T^*Q)$, let $\FDiff(\gamma)$ be the
differential operator
\begin{equation}
    \label{eq:FDiff}
    (\FDiff(\gamma)f)(\alpha_q) = \frac{d}{dt} f(\alpha_q + \gamma(q))
    \big|_{t=0}
\end{equation}
of fiber differentiation along $\gamma$, where 
$f \in C^\infty(T^*Q)$, $\alpha_q \in T^*_qQ$, and $q \in Q$. Since
all the $\FDiff(\gamma)$ commute, $\FDiff$ can be extended uniquely to
an injective algebra homomorphism from 
$\Gamma^\infty(\bigvee^\bullet T^*Q)$ into the algebra of differential
operators of $C^\infty(T^*Q)$, where zero forms
$u \in C^\infty(Q)$ act by multiplication by $\pi^*u$. 

Let  $\nabla$ be a torsion-free connection
on $Q$, and let $\mu \in \Gamma^\infty(|\!\bigwedge^n\!|\, T^*Q)$ be
a positive volume density. 
Using $\nabla$, we define the symmetrized covariant derivative $\SymD$ 
\cite[Eq.~(1.5)]{bordemann.neumaier.pflaum.waldmann:1998:pre},
\begin{equation}
    \label{eq:SymD}
    \SymD: \Gamma^\infty(\mbox{$\bigvee$}^\bullet T^*Q) \to
    \Gamma^\infty(\mbox{$\bigvee$}^{\bullet+1} T^*Q),
\end{equation}
which is a derivation of the $\vee$-product. Finally, let $\Delta$ be the
Laplacian operator on $C^\infty(T^*Q)$ coming from the indefinite
Riemannian metric on $T^*Q$ induced by the natural pairing of
vertical and horizontal spaces with respect to $\nabla$. Locally, in a
bundle chart, we have
\begin{equation}
    \label{eq:Delta}
    \Delta = \sum_k \frac{\partial^2}{\partial p_k \partial q^k} +
    \sum_{k,l,j} p_l \; \pi^*\Gamma^l_{jk} 
    \frac{\partial^2}{\partial p_j \partial p_k} +
    \sum_{k,j} \pi^* \Gamma^j_{jk} \frac{\partial}{\partial p_k},
\end{equation}
where $\Gamma^l_{jk}$ denote the Christoffel symbol of $\nabla$. 

These operators provide a nice description of the usual (formal)
differential operator calculus on $C^\infty(Q)$ in standard and in
$\kappa$-ordering, see \cite[Sect.~6]{bordemann.neumaier.waldmann:1998} and
\cite[Sect.~2]{bordemann.neumaier.pflaum.waldmann:1998:pre}.
\begin{definition}
    \label{definition:srep}
    The standard-ordered representation of a formal symbol 
    $f \in C^\infty (T^*Q)[[\lambda]]$ acting as formal series of
    differential operators on a formal wave function 
    $u \in C^\infty (Q)[[\lambda]]$ is defined by
    \begin{equation}
        \label{eq:srep}
        \srep(f) u = \iota^* \FDiff(\exp(-\im\lambda\SymD)u) f,
    \end{equation}
    where $\iota:Q \hookrightarrow T^*Q$ is the zero-section
    embedding.
\end{definition}
\begin{lemma}
    \label{lemma:stars}
    For a choice of $\nabla$ on $Q$, the expression
        \begin{equation}
        \label{eq:stars}
        \srep(f \stars g) = \srep(f) \srep(g)
    \end{equation}
    for $f, g \in C^\infty(T^*Q)[[\lambda]]$,
    defines a differential star product on $T^*Q$ 
    of standard-order type, i.e.$(\pi^*u)\stars f=(\pi^*u)f$.
\end{lemma}
The star product $\stars$ in not Hermitian, but this can be fixed as
follows. Let $\alpha \in \Gamma^\infty(T^*Q)$ be such that 
$\nabla_X \mu = \alpha(X)\mu$ for $X \in \Gamma^\infty(TQ)$, and
consider the equivalence transformation
\begin{equation}
    \label{eq:nkappa}
    N_\kappa := \eu^{-\im\kappa\lambda(\Delta + \FDiff(\alpha))}
\end{equation}
for $\kappa \in \mathbb{R}$. 
\begin{definition}
    \label{definition:stark}
    The $\kappa$-ordered star product
    $\stark$ is defined by
    \begin{equation}
        \label{eq:stark}
        f \stark g = N_\kappa^{-1} (N_\kappa f \stars N_\kappa g),
    \end{equation}
    and the corresponding $\kappa$-ordered representation on wave
    functions is defined by
    \begin{equation}
        \label{eq:krep}
        \varrho_\kappa (f) u = \srep(N_\kappa(f)) u.
    \end{equation}
    The Weyl-ordered star product is $\starw = \star_{1/2}$, and the
    Schr\"{o}dinger representation is $\wrep = \varrho_{1/2}$. We also
    set $N = N_{1/2}$.
\end{definition}
One can check that $\starw$ is Hermitian, and
the Schr\"{o}dinger representation $\wrep$ yields a
$^*$-representation of $(C^\infty(T^*Q)[[\lambda]],\starw)$ on the pre-Hilbert space 
$C^\infty_0 (Q)[[\lambda]]$ over $\mathbb{C}[[\lambda]]$ 
(see Definition \ref{def:prehilbert}) with the
usual $L^2$-inner product induced by $\mu$.
\begin{lemma}
    \label{lemma:starkproperties}
    Let $u, v \in C^\infty(Q)[[\lambda]]$ and 
    $f \in C^\infty(T^*Q)[[\lambda]]$. Then
    \begin{equation}
        \label{eq:starkpi}
        \pi^*u \stark f = \FDiff(\exp(\im\kappa\lambda\SymD)u)f
        \quad
        \textrm{and}
        \quad
        f \stark \pi^*u = \FDiff(\exp(-\im(1-\kappa)\lambda\SymD)u)f.
    \end{equation}
    In particular, $\pi^*u \stark \pi^*v = \pi^*(uv)$ whence 
    $\Exp(\pi^*u) = \pi^* \eu^u$. 
\end{lemma}

For $A \in \Gamma^\infty (T^*Q)[[\lambda]]$, let us define the operator
\begin{equation}
    \label{eq:deltaA}
    \delta_\kappa [A] 
    = \FDiff\left(
        \frac{\eu^{\im\kappa\lambda\SymD}-\eu^{-\im(1-\kappa)\lambda\SymD}} 
        {\SymD} A 
    \right).
\end{equation}
It is simple to check, using \eqref{eq:starkpi}, that it provides a 
generalization of the $\stark$-commutator with a function $\pi^*u$, i.e.
\begin{equation}
    \label{eq:deltakdu}
    \delta_\kappa[du] = \ad\nolimits_\kappa (\pi^*u).
\end{equation}
Moreover, $A \mapsto \delta_\kappa[A]$ is linear and all
$\delta_\kappa[A]$ commute.

\subsection{Deformed vector bundles over $T^*Q$ and magnetic monopoles}
\label{sec:magnetic}

We now consider deformation quantization of vector bundles over $T^*Q$
with respect to the star products $\star_{\kappa}$. As we will see,
explicit formulas for the deformed structures are obtained in this
case. We will restrict our attention to deformations of pulled-back
vector bundles $\pi^*E\to T^*Q$, where $E\to Q$, since any vector
bundle $F\to T^*Q$ is isomorphic to one of this type. For the same
reason, we assume that the Hermitian fiber metric on $\pi^*E$ is of
the form $\pi^*h_0$, for a Hermitian fiber metric $h_0$ on $E$.

Let $\{\mathcal{O}_\alpha\}$ be a good cover of $Q$, and
$\{T^*\mathcal{O}_\alpha\}$ be the corresponding good cover of $T^*Q$. 
We fix local frames $e_\alpha = \pi^*\epsilon_\alpha$ on
$T^*\mathcal{O}_\alpha$ induced by local frames 
$\epsilon_\alpha 
= (\epsilon_{\alpha, 1}, \ldots, \epsilon_{\alpha, k})$  
of $E$ on $\mathcal{O}_\alpha$. Clearly, if
$\varphi_{\alpha\beta} \in C^\infty (\mathcal{O}_{\alpha\beta})$
are transition matrices for $E$, then
$\phi_{\alpha\beta} = \pi^*\varphi_{\alpha\beta}$ are the transition
matrices for $\pi^*E$ corresponding to the frames $e_\alpha$.
\begin{proposition}
    \label{proposition:deformVBonTQ}
    Let $E \to Q$ be a complex vector bundle and $\pi^*E \to T^*Q$ its
    pull-back to $T^*Q$. Then we have:
    \begin{enumerate}
    \item The classical transition matrices
        $\phi_{\alpha\beta} = \pi^*\varphi_{\alpha\beta}$ satisfy the
        quantum cocycle condition
        \begin{equation}
            \label{eq:qccTQ}
            \phi_{\alpha\beta} \stark \phi_{\beta\gamma} \stark
            \phi_{\gamma\alpha} = \Unit
            \quad
            \textrm{and}
            \quad
            \phi_{\alpha\beta} \stark \phi_{\beta\alpha} = \Unit.
        \end{equation}
    \item For $s \in \Gamma^\infty (\pi^*E)[[\lambda]]$ and 
        $f \in C^\infty (T^*Q)[[\lambda]]$,
        \begin{equation}
            \label{eq:bulletk}
            s \bullet_\kappa f \big|_{T^*\mathcal{O}_\alpha}
            := e_\alpha (s_\alpha \stars N_\kappa (f))
            = e_\alpha N_\kappa (N_\kappa^{-1}(s_\alpha) \stark f)
        \end{equation}
        defines a global deformation quantization $\bullet_\kappa$ of
        $\pi^*E$ with respect to $\stark$ for all $\kappa$.
    \item The quantum transition matrices $\qphi_{\alpha\beta}$ with
        respect to $\bullet_\kappa$ corresponding to the frame 
        $\qe_\alpha = e_\alpha = \pi^*\epsilon_\alpha$ are 
        $\qphi_{\alpha\beta} = \phi_{\alpha\beta} =
        \pi^*\varphi_{\alpha\beta}$, for all $\kappa$.
        The local quantum trivialization isomorphisms
        $\qpsi_\alpha^{(\kappa)}$ are given by
        \begin{equation}
            \label{eq:kappaTriv}
            \qpsi_\alpha^{(\kappa)} (s) 
            = \qs_\alpha^{(\kappa)} 
            = N_\kappa^{-1} (s_\alpha),
        \end{equation}
        where $s = e_\alpha s_\alpha$ locally.
    \end{enumerate}
\end{proposition}
\begin{proof}
    The first part is clear. For the second part, let us first
    consider standard-order. In this case, 
    $\phi_{\alpha\beta} \stars s_\beta = \phi_{\alpha\beta} s_\beta$
    by \eqref{eq:starkpi} whence \eqref{eq:bulletk} is
    well-defined for $\kappa = 0$. The general case follows from
    $s \bullet_\kappa f = s \bullets N_\kappa (f)$. A local
    computation shows that \eqref{eq:bulletk} defines a deformation
    quantization of $\pi^*E$. The third part again follows from
    \eqref{eq:starkpi} and the fact that $N_\kappa \pi^* = \pi^*$. 
    
\end{proof}

In the Weyl-ordered case $\bulletw = \bullet_{1/2}$, we can also
deform the Hermitian metric $\pi^*h_0$ of $\pi^*E$. To this end we
assume that the undeformed frames $e_\alpha = \pi^*\epsilon_\alpha$
are orthonormal with respect to $\pi^*h_0$.
\begin{lemma}
    \label{lemma:HermWeylCase}
    Let $E \to Q$ be a Hermitian vector bundle with fiber metric $h_0$,
    and consider its pull back $(\pi^*E, \pi^*h_0)$. Assume
    that $e_\alpha = \pi^*\epsilon_\alpha$ are local orthonormal
    frames, and consider the Weyl-ordered deformation quantization
    $\bulletw$ of $\pi^*E$. The following holds.
    \begin{enumerate}
    \item For $s, s' \in \Gamma^\infty(\pi^*E)[[\lambda]]$,
        \begin{equation}
            \label{eq:WeylHerm}
            \qh (s, s') \big|_{T^\mathcal{O}_\alpha} := 
            \left(\qs_\alpha^{(\TinyW)}\right)^* 
            \starw {\qs'_\alpha}^{(\TinyW)}
            = \left(N^{-1} s_\alpha\right)^* \starw N^{-1} s_\alpha'
        \end{equation}
        defines a global deformation quantization of $\pi^*h_0$ with
        respect to $\bulletw$. In particular, for pulled-back sections,
        one has
        $\qh(\pi^*\sigma, \pi^*\sigma') = \pi^* h_0 (\sigma,\sigma')$.
    \item The frames $e_\alpha = \pi^*\epsilon_\alpha$ are orthonormal
        with respect to $\qh$, and hence the transition matrices are
        unitary:
        \begin{equation}
            \label{eq:WeylUnitaryPhi}
            \phi_{\alpha\beta}^* \starw \phi_{\alpha\beta} 
            = \Unit.
        \end{equation}
    \end{enumerate}
\end{lemma}
\begin{proof}
    Since \eqref{eq:WeylUnitaryPhi} is obviously satisfied, 
     \eqref{eq:WeylHerm} is globally defined. The
    remaining properties of a deformation quantization of $\pi^*h_0$
    are easily verified from the local formula. Again 
    $N_\kappa \pi^* = \pi^*$ and \eqref{eq:starkpi} imply that $\qh$
    coincides with $\pi^*h_0$ on pulled-back sections. Thus the
    $e_\alpha$ are still orthonormal.
    
\end{proof}

Let us now consider a line bundle $L \to Q$, 
with pull-back $\pi^*L \to T^*Q$. In this case, we can describe
the deformed endomorphisms (with respect to $\bullet_\kappa$)
explicitly by using a connection $\nabla^L$ on $L$. The
frame $e_\alpha = \pi^*\epsilon_\alpha$ is a single
non-vanishing local section of $\pi^*L$, and $\nabla^L$ determines
local connection one-forms 
$A_\alpha \in \Gamma^\infty(T^*\mathcal{O}_\alpha)$ by
\begin{equation}
    \label{eq:localAalpha}
    \nabla^L_X \epsilon_\alpha = - \im A_\alpha(X) \epsilon_\alpha,
\end{equation}
where $X \in \Gamma^\infty
(TQ)$. Let $B$ be the (global) curvature two-form,
\begin{equation}
    \label{eq:Bfield}
    B = dA_\alpha.
\end{equation}
We assume $\nabla^L$ to be compatible with
$h_0$, so that the forms $A_\alpha$ and $B$ are real. Using these local
one-forms we can define local series of differential operators
$S_\alpha^{(\kappa)}$ by
\begin{equation}
    \label{eq:Salphadef}
    S_\alpha^{(\kappa)} (f) = \eu^{\im\delta_\kappa[A_\alpha]} (f).
\end{equation}
Note that the operator $S_\alpha^{(\kappa)}$ is just the
$\kappa$-ordered quantized fiber translation by the one-form $\lambda
A_\alpha$ in the sense of
\cite[Thm.~3.4]{bordemann.neumaier.pflaum.waldmann:1998:pre}.
\begin{lemma}
    \label{lemma:Skappa}
    For $\qphi_{\alpha\beta} = \phi_{\alpha\beta} =
    \pi^*\varphi_{\alpha\beta}$ 
    the relation
    \begin{equation}
        \label{eq:PhifPhi}
        \qphi_{\alpha\beta} \stark f \stark \qphi_{\beta\alpha}
        =
        \eu^{\im\delta_\kappa[A_\alpha]} \, 
        \eu^{-\im\delta_\kappa[A_\beta]} \, (f)
    \end{equation}
    holds for all $f \in C^\infty(T^*Q)[[\lambda]]$.
 \end{lemma}   
  \begin{proof}
    Choose local functions 
    $c_{\alpha\beta} \in C^\infty(\mathcal{O}_{\alpha\beta})$ such
    that $\varphi_{\alpha\beta} = \eu^{2\pi\im c_{\alpha\beta}}$. Then
    we know that $A_\alpha - A_\beta = 2\pi dc_{\alpha\beta}$ and
    \eqref{eq:PhifPhi} is a simple computation using
    \eqref{eq:deltakdu}, Lem.~\ref{lemma:starkproperties} and the
    commutativity of all $\delta_\kappa[\cdot]$. 
    
\end{proof}

As a result, (\ref{eq:ChangeQuantACoeff}) is satisfied, and hence
\begin{equation}
    \label{eq:Sfkappa}
    S_\kappa (f) s \big|_{T^*\mathcal{O}_\alpha} :=
    e_\alpha \bullet_\kappa 
    \left(
        S_\alpha^{(\kappa)}(f) \stark \qs^{(\kappa)}_\alpha
    \right) 
\end{equation}
defines a global endomorphism $S_\kappa(f)$ of 
$(\Gamma^\infty(\pi^*L)[[\lambda]],\bullet_{\kappa})$ for any 
$f \in C^\infty(T^*Q)[[\lambda]]$. Also observe that 
$S_\alpha^{(\kappa)} (\pi^*u) = \pi^*u$.

Let  $\stark'$ be the star product induced by the operator
product of deformed endomorphisms,
\begin{equation}
    \label{eq:starkB}
    f \stark' g = S_\kappa^{-1} (S_\kappa(f)S_\kappa(g))
    = \left(S_\alpha^{(\kappa)}\right)^{-1} 
    \left(S_\alpha^{(\kappa)}(f) \stark S_\alpha^{(\kappa)}(g)\right).
\end{equation}
It follows from the explicit
form  of the local equivalence transformations \eqref{eq:Salphadef} and
\cite[Thm.~4.1 and Thm.~4.6]{bordemann.neumaier.pflaum.waldmann:1998:pre}
that the star product $\stark'$ coincides with the one 
constructed in \cite{bordemann.neumaier.pflaum.waldmann:1998:pre}:
\begin{proposition}
    \label{proposition:kappaClassB}
    The star product $\stark'$ coincides with
    $\star_\kappa^{-\lambda B}$ from
    \cite[Thm.~4.1]{bordemann.neumaier.pflaum.waldmann:1998:pre}. Its
    characteristic class is given by
    \begin{equation}
        \label{eq:kappaClass}
        c\left(\star_\kappa^{-\lambda B}\right) =
        \im [\pi^* B] = 2 \pi \im \, c_1 (\pi^*L).
    \end{equation}
\end{proposition}
Note that \eqref{eq:kappaClass} is consistent with
(\ref{eq:explicitPhi}) since the characteristic class of $\stark$
vanishes, see
\cite[Thm.~4.6]{bordemann.neumaier.pflaum.waldmann:1998:pre}.
\begin{remark}
    More generally 
    \cite{bordemann.neumaier.pflaum.waldmann:1998:pre}, 
    one can explicitly construct a star product $\star^B_\kappa$,
    for any formal series of closed two-forms 
    $B \in \Gamma^\infty (\bigwedge^2 T^*Q)[[\lambda]]$, with
    $c(\star^B) = \frac{1}{\im\lambda}[\pi^*B]$. In particular, any star
    product on $T^*Q$ is equivalent to some $\star_\kappa^{B}$.
\end{remark}

The physical interpretation of the star products 
$\star^{-\lambda B}_\kappa$ is discussed in
\cite{bordemann.neumaier.pflaum.waldmann:1998:pre,waldmann:1999}:
they correspond to the quantization of a charged particle, with
electric charge $1$, moving in $Q$ under the influence of a magnetic
field $B$. With this in mind, we can think of  non-trivial
characteristic classes of star products on $T^*Q$ as corresponding to
 topologically non-trivial magnetic fields, i.e. to
the presence of magnetic monopoles.
The integral
$m = \frac{1}{4\pi} \int_{S^2} B$ gives the amount of
`magnetic charge' inside this $2$-sphere $S^2$. Thus the integrality of $B$
implies that $2m \in \mathbb{Z}$, which is  Dirac's
integrality/quantization condition for magnetic charges $m$. We summarize the
discussion:
\begin{theorem}
    \label{theorem:METQ}
    Let $B \in \Gamma^\infty(\bigwedge^2 T^*Q)[[\lambda]]$ be a
    sequence of closed two-forms, and $\star^{-\lambda B}$ the
    star product in \cite{bordemann.neumaier.pflaum.waldmann:1998:pre}.
    Then $\star_\kappa^{-\lambda B}$
    is Morita equivalent to $\stark$ if and only if
    $\frac{1}{2\pi} B$ is an integral two-form.
    In physical terms, the quantization with magnetic field $B$ is Morita
    equivalent to the quantization without magnetic field if and
    only if Dirac's integrality condition for the magnetic charge
    of $B$ is fulfilled.
\end{theorem}

This theorem suggests the physical interpretation of  characteristic
classes of star products on arbitrary symplectic manifolds as
`intrinsic magnetic monopole fields', and of Morita equivalence as
Dirac's integrality condition for the `relative fields'.

\subsection{Rieffel induction of the  Sch\"{o}dinger representation}
\label{sec:rieffel}

Let $\starw$ be the Weyl-ordered star product on $T^*Q$, and let
$\wrep$ be the Schr\"{o}dinger representation \eqref{eq:krep} of
$\starw$ on (formal) wave functions 
$\mathfrak{H} = C^\infty_0 (Q)[[\lambda]]$, 
with $L^2$-inner product coming from $\mu$, see
\cite{bordemann.neumaier.waldmann:1998,bordemann.neumaier.waldmann:1999,bordemann.neumaier.pflaum.waldmann:1998:pre}.
We now illustrate the consequences of Morita equivalence by
constructing the $^*$-representation of $\star_\TinyW^{-\lambda B}$
induced (in the sense of Rieffel induction) by $\wrep$.

Let $L \to Q$ be a Hermitian line bundle, and let $\pi^*L \to T^*Q$ be
its pull-back, endowed with a quantization $\bulletw$
and $\qh$ as before. By fixing a compatible connection
$\nabla^L$, we obtain a star product
$\star_\TinyW^{-\lambda B}$ by \eqref{eq:starkB}
such that $\Gamma^\infty (\pi^*L)[[\lambda]]$ has a bimodule structure
with respect to $\star_\TinyW^{-\lambda B}$ and $\star_\TinyW$. 
As shown in
\cite[Sect.~8 and 9]{bordemann.neumaier.pflaum.waldmann:1998:pre},
this data determines a $^*$-representation $\eta_\TinyW$ of 
$\star_\TinyW^{-\lambda B}$ on $\Gamma^\infty_0 (L)[[\lambda]]$, with
$L^2$-inner product defined by $h_0$ and the volume density $\mu$. We
have the following explicit local formula
\begin{equation}
    \label{eq:LocalRepB}
    \eta_\TinyW (f) (\epsilon_\alpha \sigma_\alpha) = 
    \epsilon_\alpha 
    \wrep \left(\eu^{\im\delta_\TinyW[A_\alpha]} f\right)\sigma_\alpha,
\end{equation}
where 
$\sigma = \epsilon_a \sigma_\alpha \in \Gamma^\infty(L)[[\lambda]]$, 
see \cite[Eq.~(5.4) and Thm.~8.2]
{bordemann.neumaier.pflaum.waldmann:1998:pre} (The missing
minus sign comes from a different convention for the Chern class of
$L$.) We shall now show that $\eta_\TinyW$ is canonically unitarily
equivalent to the Rieffel induction of the Schr\"{o}dinger
representation $\wrep$ of $\starw$.
\begin{theorem}
    \label{theorem:RieffelSchroedinger}
    Let $(\mathfrak{K}, \rho)$ be the $^*$-representation of
    $\star_\TinyW^{-\lambda B}$ obtained by Rieffel induction of the
    Schr\"{o}dinger representation $(\mathfrak{H}, \wrep)$, using the
    equivalence bimodule 
    $\mathcal{L} = \Gamma^\infty (\pi^*L)[[\lambda]]$. The following holds:
    \begin{enumerate}
    \item Let $s \in \Gamma^\infty(\pi^*L)[[\lambda]]$ and 
        $u \in C^\infty_0(Q)[[\lambda]]$. Then
        \begin{equation}
            \label{eq:Utilde}
            \widetilde{\mathfrak{K}} 
            = \mathcal{L}\otimes_{\starw} \mathfrak{H} \ni \; 
            s \otimes u \; \mapsto \; 
            \epsilon_\alpha \wrep\left(\qs_\alpha^{(\TinyW)}\right) u
            \; \in \Gamma^\infty_0 (L)[[\lambda]]            
        \end{equation}
        extends to a well-defined global $\mathbb{C}[[\lambda]]$-linear
        map $\widetilde{U}$, which is isometric and surjective.
    \item $\widetilde{U}$ induces a unitary map 
        $U: \mathfrak{K} \to \Gamma^\infty_0 (L)[[\lambda]]$.
    \item $U$ is an intertwiner between $\rho$ and $\eta_\TinyW$.
    \end{enumerate}
\end{theorem}
\begin{proof}
    Let $s = e_\alpha \bulletw \qs^{(\TinyW)}_\alpha$. A
    straightforward computation shows that
    $\epsilon_\alpha \wrep(\qs^{(\TinyW)}_\alpha) u
    = \epsilon_\beta \wrep(\qs^{(\TinyW)}_\beta) u$,
    since $\phi_{\alpha\beta} = \pi^*\varphi_{\alpha\beta}$ and
    $\wrep$ is a representation satisfying $\wrep(\pi^*v) = v$.
    Thus the right hand side of \eqref{eq:Utilde} is a global
    section. A similar computation shows that
    $\widetilde{U}(s \bulletw f \otimes u) 
    = \widetilde{U}(s \otimes \wrep(f)u)$, 
    whence $\widetilde{U}$ is well-defined. From the fact that $\wrep$
    is a $^*$-representation, one obtains for sections/functions with
    small enough support the relation
    \begin{equation}
        \label{eq:Isometric}
        \int_Q h_0 \left(
            \epsilon_\alpha \wrep\left(\qs^{(\TinyW)}_\alpha\right) u,
            \epsilon_\alpha \wrep\left(\qt^{(\TinyW)}_\alpha\right) v
        \right) \; \mu
        =
        \int_Q \cc{u} \, \wrep(\qh(s,t)) v \; \mu.
    \end{equation}
    Then a partition of unity argument implies that $\widetilde{U}$ is
    isometric. Finally we choose $\sigma \in \Gamma^\infty_0 (L)$ and
    $u$ such that $u = 1$ on $\supp \sigma$. Then clearly
    $\widetilde{U}(\pi^*\sigma \otimes u) = \sigma u = \sigma$ implies
    surjectivity. This shows the first part. The second part is trivial
    since $\mathfrak{K}$ is the quotient of $\widetilde{\mathfrak{K}}$
    by the vectors of length zero. For the third part, we compute
    locally 
    \begin{equation*}
        \label{eq:UIntertwiner}
        \begin{split}
            & U(\rho(f) s \otimes u) \\
            & = \epsilon_\alpha \wrep(S_\alpha^{(\TinyW)}(f) \starw
            \qs^{(\TinyW)}_\alpha) u
            = \epsilon_\alpha 
            \wrep\left(\eu^{\im\delta_\TinyW[A_\alpha]}(f)\right)
            \wrep\left(\qs^{(\TinyW)}_\alpha\right) u
            = \eta_\TinyW (f) U (s \otimes u),       
        \end{split}
    \end{equation*}
    which is sufficient since all representations are local.
    
\end{proof}

The $^*$-representation $\eta_\TinyW$ is well-known, for instance,
from geometric quantization \cite[Sect.~8.4]{woodhouse:1992}: It is
precisely the representation obtained if the symplectic form satisfies
the integrality condition of pre-quantization. The difference is that
we have treated $\hbar$ as a formal parameter $\lambda$, so the
correction to the canonical symplectic form occurs in first order of
$\lambda$. For a further discussion see also
\cite{bordemann.neumaier.pflaum.waldmann:1998:pre}.

As we just saw, $\eta_\TinyW$ can be obtained as a result of Rieffel
induction applied to the ordinary Schr\"{o}dinger representation
$\wrep$. We remark that, by Morita equivalence, $\starw$ and
$\star_\TinyW^{-\lambda B}$ have equivalent categories of
$^*$-representations, and the correspondence of $\wrep$ and
$\eta_\TinyW$ is just one example of this more general fact. These
considerations are based on the approach to quantization where primary
objects are observable algebras, as opposed to specific
$^*$-representations.

The results in this paper illustrate that several constructions and
techniques present in more analytic approaches to quantization find
counterparts in formal deformation quantization. It is interesting to
investigate how far one can go without convergence.

%
%

\appendix

\section{Star exponentials and star logarithms}
\label{sec:explog}

In this appendix we recall a few properties of the star
exponential \cite{bayen.et.al:1978} and the star logarithm (see
\cite{bordemann.roemer.waldmann:1998,waldmann:1999} for details).

Let $\star$ be a star product on a Poisson manifold $M$. Let
$H = \sum_{r=0}^\infty \lambda^r H_r \in C^\infty(M)[[\lambda]]$ and
consider the differential equation
\begin{equation}
    \label{eq:ExpDiff}
    \frac{d}{dt} f(t) = H \star f(t),
    \qquad
    f(0) = 1,
\end{equation}
for $t \in \mathbb{R}$ and  $f(t) \in C^\infty(M)[[\lambda]]$. The
next result follows from 
\cite[Lem.~2.2, 2.3]{bordemann.roemer.waldmann:1998} and
\cite[Sect.~1.4.2]{waldmann:1999}:
\begin{proposition}
    \label{proposition:Exp}
    For any $H \in C^\infty(M)[[\lambda]]$, the differential 
    equation \eqref{eq:ExpDiff} has a unique solution, denoted by
    $t \mapsto \Exp(tH)$, satisfying the following properties:
    \begin{enumerate}
    \item $\Exp(tH) = \sum_{r=0}^\infty \lambda^r \Exp(tH)_r$, with
        $\Exp(tH)_0 = \eu^{tH_0}$ and $\Exp(H)_{r+1}$ equals 
        $\eu^{H_0} H_{r+1}$ plus terms only depending on 
        $H_0$, \ldots $H_r$.
    \item $\Exp(tH) \star H = H \star \Exp(tH)$, and
        $\Exp((t+t')H) = \Exp(tH) \star \Exp(t'H)$.
    \item If $\star$ is a Hermitian star product, then 
        $\cc{\Exp(tH)} = \Exp(t\cc{H})$.
    \item $\Exp(H) = 1$ if and only if $H$ is constant on each
        connected component of $M$ and equal to $2\pi\im k$ for some
        $k \in \mathbb{Z}$.
    \item For all $f \in C^\infty(M)[[\lambda]]$ we have
        \begin{equation}
            \label{eq:ExpAd}
            \eu^{\ad(H)} (f) = \Exp(H) \star f \star \Exp(-H),
        \end{equation}
        where $\ad(H) = [H, \cdot]_\star$ denotes the
        $\star$-commutator.
    \end{enumerate}
\end{proposition}

For $f, g \in C^\infty(M)[[\lambda]]$, consider the Baker-Campbell-Hausdorff formula
\begin{equation}
    \label{eq:BCH}
    f \circ_\star g 
    = f + g + \frac{1}{2} [f, g] 
    + \frac{1}{12}\left([f, [f,g]] + [g, [g,f]]\right) + \cdots
\end{equation}
Since in zeroth order
the star commutator vanishes, the series \eqref{eq:BCH} is a
well-defined formal power series in
$C^\infty(M)[[\lambda]]$, and one has
\begin{equation}
    \label{eq:ExpBCH}
    \Exp(f) \star \Exp(g) = \Exp(f \circ_\star g).
\end{equation}
See e.g. \cite[Lem.~4.1]{gutt.rawnsley:1999} for the properties of $\circ_\star$.

More generally, we define star logarithms in the following
way. Let $U \subseteq M$ be a contractible open subset, and let 
$f = \sum_{r=0}^\infty \lambda^r f_r \in C^\infty(U)[[\lambda]]$. If
$f_0(x) \ne 0$ for all $x \in U$, then there exists a smooth logarithm
$H_0 = \ln(f_0) \in C^\infty(U)$ for the pointwise product, 
unique up to constants in 
$2\pi\im \mathbb{Z}$. If we have fix the choice of the
classical $\ln$, then Prop.~\ref{proposition:Exp} ensures that we can
find $H_1, H_2, \ldots \in C^\infty(U)$ by recursion such that
$\Exp(H) = f$ for $H = \sum_{r=0}^\infty \lambda^r H_r$. We
write $H = \Ln(f)$, and call it (the/a) star
logarithm of $f$ corresponding to the choice of the classical
$\ln(f_0)$. Again $H$ is unique up to constants in $2\pi\im\mathbb{Z}$
and
\begin{equation}
    \label{eq:ExpLogLogExp}
    \Exp(\Ln(f)) = f
    \quad
    \textrm{and}
    \quad
    H = \Ln(\Exp(H)) \;\textrm{mod}\; 2\pi\im\mathbb{Z}.
\end{equation}

\section{Rieffel induction and strong Morita equivalence}
\label{sec:strong}

This appendix recalls the notions of algebraic Rieffel induction and
strong Morita equivalence for $^*$-algebras over an ordered ring.
For simplicity, we assume $^*$-algebras to be unital.
The reader is referred to
\cite{bursztyn.waldmann:2001a,bursztyn.waldmann:2000a:pre} for
details. 

Let $\ring R$ be an ordered ring, and let 
$\ring C = \ring R(\im)$ with $\im^2 = -1$. 
The main examples from deformation quantization are 
$\ring R=\mathbb R$ and $\ring R = \mathbb{R}[[\lambda]]$. We consider
the following generalization of complex pre-Hilbert spaces.
\begin{definition}\label{def:prehilbert}
    Let $\mathfrak{H}$ be a $\ring C$-module. A 
    \emph{Hermitian inner product} on $\mathfrak H$ is a sesquilinear
    map 
    $\SP{\cdot,\cdot}: \mathfrak{H} \times \mathfrak{H} \to \ring C$
    such that $\SP{\phi,\psi} = \cc{\SP{\psi,\phi}}$, and
    $\SP{\phi,\phi} > 0$ for all $\phi \ne 0$. The pair $(\mathfrak H,
    \SP{\cdot,\cdot})$ is called a 
    \emph{pre-Hilbert space over $\ring C$}.
\end{definition}
Let $\mathfrak{B}(\mathfrak{H})$ be the $^*$-algebra of adjointable
$\ring C$-linear endomorphisms of $\mathfrak{H}$. If $\mathcal A$ is a
$^*$-algebra over $\ring C$, a \emph{$^*$-representation} of 
$\mathcal A$ on $\mathfrak H$ is a $^*$-homomorphism 
$\pi: \mathcal{A} \to \mathfrak{B}(\mathfrak{H})$. We denote the
category of nondegenerate (i.e. $\pi(\Unit) = \id$)
$^*$-representations of $\mathcal A$ by $\sRep(\mathcal{A})$. 
Following the analogy with $C^*$-algebras, we consider
\cite{bordemann.waldmann:1998}: 
\begin{definition}
    A $\ring C$-linear functional 
    $\omega: \mathcal{A} \to \ring C$ is called \emph{positive} if
    $\omega(A^*A) \ge 0$ for all $A \in \mathcal{A}$. An element
    $A \in \mathcal A$ is called \emph{positive} if 
    $\omega(A) \ge 0$ for all positive linear functionals $\omega$.
\end{definition}
Elements of the form 
$A = b_1 B_1^* B_1 + \cdots + b_n B_n^*B_n$, with $b_i > 0$ and 
$B_i \in \mathcal{A}$, are necessarily positive, and called
\emph{algebraically positive}. These definitions recover the usual
notions of positivity on $C^*$-algebras. If 
$\mathcal A = C^\infty(M)$, then positive linear functionals
are compactly supported positive measures, and positive elements are
usual positive functions.

To describe Rieffel induction \cite{rieffel:1974}, 
we consider algebraic analogs of Hilbert modules.
\begin{definition}
    Let $\mathcal{E}$ be a $\mathcal{A}$-right module.
    An \emph{$\mathcal{A}$-valued Hermitian inner product} is a 
    $\ring C$-sesquilinear map 
    $\SP{\cdot,\cdot}: \mathcal{E} \times \mathcal{E} \to \mathcal{A}$
    such that $\SP{x,y} = \SP{y,x}^*$, $\SP{x,y\cdot A} = \SP{x,y}A$,
    and $\SP{x,x}$ is positive in $\mathcal{A}$. 
\end{definition}
Suppose $\mathcal E$ is a ($\mathcal{B}$,$\mathcal{A}$)-bimodule, 
equipped with an $\mathcal{A}$-valued Hermitian inner product, so that
\begin{equation}
    \label{eq:BByAdjoints}
    \SP{B \cdot x, y} = \SP{x, B^* \cdot y}.
\end{equation}
Let $(\mathfrak{H}, \pi)$ be a $^*$-representation of $\mathcal{A}$. 
Consider the space $\widetilde{\mathfrak{K}} 
= \mathcal{E} \otimes_{\mathcal{A}} \mathfrak{H}$,
endowed with its canonical $\mathcal{B}$-left module
structure, and set
\begin{equation}
    \label{eq:RieffelInnerProd}
    \SP{x \otimes \phi, y \otimes \psi} 
    =
    \SP{\phi, \pi(\SP{x,y}) \psi}.
\end{equation}
We assume that $\mathcal E$ is such that (\ref{eq:RieffelInnerProd})
defines a positive semi-definite inner product on
$\widetilde{\mathfrak{K}}$ for all $^*$-representations
(this is always the case for $C^*$-algebras and
for star product algebras if $\mathcal{E}$ is a deformation
quantization of a Hermitian vector bundle
\cite{bursztyn.waldmann:2000a:pre}). Factoring 
$\widetilde{\mathfrak{K}}$ by the vectors of length zero, we obtain a
pre-Hilbert space $\mathfrak{K}$ over $\ring C$ equipped with a
$^*$-representation of $\mathcal{B}$. This induced $^*$-representation
is denoted by $\mathsf{R}_{\mathcal{E}}(\mathfrak H,\pi)$, and the
induction process is functorial. 
\begin{definition}
    \label{definition:rieffel}
    The functor $\mathsf{R}_{\mathcal{E}}: 
    \sRep(\mathcal{A}) \to \sRep(\mathcal{B})$
    is called Rieffel induction.
\end{definition}
In order to get an equivalence of categories, we assume that 
$\mathcal E$ is, in addition, equipped with a $\mathcal{B}$-valued
Hermitian inner product
$\Theta_{\cdot,\cdot}: \mathcal{E}\times\mathcal{E} \to  \mathcal{B}$
so that $\Theta_{x,y\cdot A} = \Theta_{x\cdot A^*, y}$. We require the
compatibility 
\begin{equation}
    \label{eq:ThetaSPComp}
    \Theta_{x,y} \cdot z = x \cdot \SP{y,z},
\end{equation}
and assume that the following \emph{fullness conditions} hold:
\begin{equation}
    \label{eq:Fullness}
    \begin{array}{rcl}
        \mathcal{A} 
        & = & \ring C \textrm{-span} 
        \{ \SP{x,y} \; | \; x, y \in \mathcal{E} \}
        \\
        \mathcal{B}
        & = & \ring C \textrm{-span}
        \{ \Theta_{x,y} \; | \; x, y \in \mathcal{E} \}.
    \end{array}
\end{equation}
\begin{definition}
    \label{definition:moritabimod}
    A ($\mathcal{B}$,$\mathcal{A}$)-bimodule $\mathcal E$ equipped with
    full $\mathcal{A}$- and $\mathcal{B}$-valued inner products
    satisfying the above properties is called an
    equivalence bimodule, and the $^*$-algebras $\mathcal{A}$
    and $\mathcal{B}$ are called strongly Morita equivalent as
    $^*$-algebras over $\ring C$. 
\end{definition}
\begin{proposition}
    \label{lemma:SME}
    Let $\mathcal{A}$, $\mathcal{B}$ be strongly Morita equivalent
    unital $^*$-algebras over
    $\ring C$, with equivalence bimodule $\mathcal E$. Then
    $\mathsf{R}_{\mathcal{E}}: 
    \sRep(\mathcal{A}) \to \sRep(\mathcal{B})$
    is an equivalence of categories.
\end{proposition}
\begin{remark}
    \begin{enumerate}
    \item The bimodule $\mathcal E$
        is also an equivalence bimodule in the purely ring theoretic sense
        of Morita equivalence \cite[Sec.~7]{bursztyn.waldmann:2001a}. In
        particular, $\mathcal{E}$ is finitely generated and projective over
        $\mathcal A$ and $\mathcal B$.
    \item Analogous results hold for nonunital $^*$-algebras. In particular,
        if $\mathcal A$ and $\mathcal B$ are $C^*$-algebras, then they are
        strongly Morita equivalent (in the usual sense of operator
        algebras \cite{rieffel:1982}) if and only if their Pedersen ideals are
        strongly Morita equivalent in the sense of
        Definition~\ref{definition:moritabimod}
        \cite[Sec.~3]{bursztyn.waldmann:2000a:pre}. 
    \end{enumerate}
\end{remark}

%
%

\begin{footnotesize}

\end{footnotesize}

\end{document}